\documentclass[11pt,twoside,letter]{article}

\usepackage{hyperref}

\usepackage{fancyhdr}

\usepackage{graphicx}

\usepackage{amsmath,amssymb,amscd,amsthm}

\def\l{\lambda}

\def\e{\epsilon}
\def\o{\omega}

\def\bp{\begin{proposition}}
\def\ep{\end{proposition}}
\def\bt{\begin{theo}}
\def\et{\end{theo}}
\def\be{\begin{equation}}
\def\ee{\end{equation}}
\def\bl{\begin{lemma}}
\def\el{\end{lemma}}
\def\bc{\begin{corollary}}
\def\ec{\end{corollary}}
\def\pr{{\bf Proof: }}

\def\bd{\begin{definition}}
\def\ed{\end{definition}}
\def\C{{\mathbb C}}

\def\arg{{\rm arg\,}}
\def\min{{\rm min\,}}
\def\max{{\rm max\,}}

\def\arg{{\rm arg\,}}
\def\min{{\rm min\,}}
\def\max{{\rm max\,}}
\def\l{\ell}

\newtheorem{theo}{Theorem}[section]
\newtheorem{lemma}{Lemma}[section]
\newtheorem{corollary}{Corollary}[section]
\newtheorem{proposition}{Proposition}[section]

\theoremstyle{definition}
\newtheorem{definition}{Definition}[section]

\newtheorem{remark}{Remark}[section]

\renewcommand{\vec}[1]{\boldsymbol{#1}}

\begin{document}

\title{Accuracy of Algebraic Fourier Reconstruction for Shifts of Several Signals}

\author{Dmitry Batenkov$^{\dagger\; \sharp}$, Niv Sarig$^{\dagger}$, Yosef Yomdin$^{\dagger}$  \\
\small \{dima.batenkov,niv.sarig,yosef.yomdin\}@weizmann.ac.il}
%
%
\date{}

%
%
\maketitle
\thispagestyle{fancy}

%
%
\markboth{\footnotesize \rm \hfill D. BATENKOV AND N. SARIG AND Y. YOMDIN \hfill}
{\footnotesize \rm \hfill Accuracy of Algebraic Fourier Reconstruction for Shifts of Several Signals \hfill}

\thanks{\small This research was supported by the Adams Fellowship Program of the Israeli Academy of Sciences and Humanities,
ISF Grant No. 779/13, and by the Minerva foundation.}

{$^{\dagger}$ \small Department of Mathematics, Weizmann Institute of Science, Rehovot 76100, Israel \\
$^{\sharp}$ \small Department of Mathematics, Ben-Gurion University of the Negev, Beer-Sheva 84105, Israel}

%
%

\begin{abstract}
We consider the problem of ``algebraic reconstruction'' of linear combinations of shifts of several known signals
$f_1,\ldots,f_k$ from the Fourier samples. Following \cite{Bat.Sar.Yom2}, for each $j=1,\ldots,k$ we choose sampling set
$S_j$ to be a subset of the common set of zeroes of the Fourier transforms ${\cal F}(f_\l), \ \l \ne j$, on which
${\cal F}(f_j)\ne 0$. It was shown in \cite{Bat.Sar.Yom2} that in this way the reconstruction system is ``decoupled'' into $k$
separate systems, each including only one of the signals $f_j$. The resulting systems are of a ``generalized Prony'' form.

However, the sampling sets as above may be non-uniform/not ``dense enough'' to allow for a unique reconstruction of the shifts 
and amplitudes. In the present paper we study uniqueness and robustness of non-uniform Fourier sampling of signals as above,
investigating sampling of exponential polynomials with purely imaginary exponents. As the main tool we apply a well-known 
result in Harmonic Analysis: the Tur\'an-Nazarov inequality (\cite{Naz}), and its generalization to discrete sets, obtained 
in \cite{Fri.Yom}. We illustrate our general approach with examples, and provide some simulation results.
\vspace{5mm} \\
 {\it Key words and phrases} : non-uniform sampling, Tur\'an-Nazarov inequality, exponential fitting, Prony systems
\vspace{3mm}\\
 {\it 2010 AMS Mathematics Subject Classification} --- 94A20, 65T40  
\end{abstract}


\section{Introduction}
\setcounter{equation}{0}

In this paper we investigate robustness of Fourier reconstruction of signals of the following a priori known form:

\be\label{equation_decoupling_model}
F(x)=\sum_{j=1}^k \sum_{q=1}^{q_j} a_{jq} f_j(x-x_{jq}),
\ee
with $a_{jq}\in \mathbb{C}, \ x_{jq} \in {\mathbb R}.$ We assume that the signals $f_1,\dots, f_k: \mathbb{R}\to\mathbb{R}$ are known
(in particular, their Fourier transforms ${\cal F}(f_j)$ are known), while $a_{jq}, \ x_{jq}$ are the unknown signal parameters, which
we want to find from Fourier samples of $F$.

Practical importance of signals of the form \eqref{equation_decoupling_model} is well recognized in the literature. For instance, they appear in digital processing of neuronal signals, bioimaging, image processing and ultrawideband communications \cite{Vet,gedalyahu2011multichannel}. They are of relevance also in inverse moment problems, an important subject in mathematical physics \cite{gust.Moments, peter2011nonlinear}. 

We follow a general line of the ``Algebraic Sampling'' approach (see \cite{Bat.Yom1,sig.ack,Vet} and references therein), i.e. we
reconstruct the values of the unknown parameters, solving a system of non-linear equations, imposed by the measurements. The equations in
this system appear as we equate the ``symbolic'' expressions of the Fourier samples, obtained from (\ref{equation_decoupling_model}), to
their actual measured values.

Our specific strategy, as suggested in \cite{Bat.Sar.Yom2,Sar}, is as follows: we choose a sampling set
$S_j \subset {\mathbb R}, \  j=1,\ldots,k,$ in a special way, in order to ``decouple'' the reconstruction system, and to reduce
it to $k$ separate systems, each including only one of the signals $f_j$. To achieve this goal we take $S_j$ to be a subset of the sets $W_j$ 
of common zeroes of the Fourier transforms ${\cal F}(f_\l), \ \l\ne j$. It was shown in \cite{Bat.Sar.Yom2} that the decoupled systems turn
out to be exactly the same as those which appear in the fitting of exponential polynomials on sets $S_j$ (systems (\ref{equation_decoupled1}) 
in Section \ref{Four.Main} below).

In this paper we restrict ourselves to one-dimensional case. A presentation of the Fourier
Decoupling method in several variables, as well as some initial uniqueness results, can be found in \cite{Sar,Bat.Sar.Yom2}. On the other
hand, we explicitly assume here that $k\geq 2$. So the usual methods which allow one to solve this problem ``in closed form'' in
the case of shifts of a single function (see \cite{Vet,Bat.Sar.Yom1,Sar}) are not directly applicable. Still, as it was shown in
\cite{Bat.Sar.Yom2}, in many cases an explicit reconstruction from a relatively small collection of Fourier samples of $F$ is possible. Let us also stress that the decoupling method of \cite{Sar,Bat.Sar.Yom2} in dimension one can ``generically'' be applied only to the shifts of at most 
two different signals. Indeed, for three or more signals the sampling sets $S_j$ are the intersections of at least two different discrete sets
(the sets of zeroes of the Fourier transforms ${\cal F}(f_\l), \ \l\ne j$), so ``generically'' $S_j$ are empty. However, in many important
``non-generic'' situations of $k>2$ one-dimensional signals the resulting sampling sets are dense enough for a robust reconstruction. Accordingly,
our main result - Theorem \ref{sign.main} below - is stated for an arbitrary $k$.

If the points $s_{j\l}\in S_j, \quad \l=1,2,\dots,$ form an arithmetic progression, the reconstruction systems (\ref{equation_decoupled1}) are
very closely related to the standard Prony system (see, for instance, \cite{Bat.Yom} and discussion therein). However, the sampling sets $S_j$,
being subsets of the sets $W_j$ of zeroes of the Fourier transforms ${\cal F}(f_\l), \ \l\ne j$, are completely defined by the original signals
$f_\l$, and cannot be altered in order to make sampling more stable. These sets usually are non-uniform, therefore the standard methods for 
robust solution of Prony systems cannot be applied. Even if $S_j$ forms an arithmetic progression, it may turn out to be ``insufficiently dense'' 
to allow a robust reconstruction of the shifts and amplitudes (see an example in Section \ref{examples} below). Because of these reasons, we 
restrict ourselves to only one solution method for system (\ref{equation_decoupled1}) - that of the least squares fitting, mainly because of its 
relative insensitivity to the specific geometry of the sampling set. Accordingly, we do not consider in this paper other approaches, which can be 
more efficient in certain specific circumstances. Let us only mention that non-uniform sampling is an active area of research, see e.g. \cite{adk,Marvasti} and references therein.

The main goal of the present paper is to study uniqueness and robustness of the Fourier decoupling method. We define a ``metric span'' $\omega(S)$ 
of sampling sets $S$, which is a simple geometric quantity, taking into account both the geometry of $S$, as well as the maximal shifts allowed in
the signal $F$ (which are the maximal frequencies of the exponential polynomials appearing in the Fourier transform of $F$). Our main results - 
Theorem \ref{sign.main} and Corollary \ref{density} below - provide, in terms of the metric span $\omega$ a ``density-like'' geometric condition 
on the common sets $W_j$ of zeroes of the Fourier transforms ${\cal F}(f_\l), \ \l\ne j$, which, in the case of no noise, guarantees uniqueness of 
the least square reconstruction via the decoupled systems. In the noisy case Theorem \ref{sign.main} provides an upper bound for the maximal error 
of the least square reconstruction. The  proof of these results relies on a stability estimate for non-uniform sampling of exponential polynomials, whose derivation constitutes the bulk of the paper (Section \ref{Exp.Pol}). The principal result there, which might be of independent interest, is Theorem \ref{non.unif.rec}, whose proof is, in turn, based on the classical Tur\'an-Nazarov 
inequality \cite{Naz}, and its generalization to discrete sets, obtained recently in \cite{Fri.Yom}.

We hope that our results will provide  useful criteria for applicability of Fourier Decoupling method to  specific models of the form \eqref{equation_decoupling_model} in specific applications, and a guiding principle  for designing relevant sampling strategies. While the stability bounds in Theorem \ref{sign.main}  increase exponentially in the number of shifts $q_j$, this seems to be consistent with the general principle that reconstruction methods based on sparsity are poorly conditioned with respect to model complexity, see \cite{donoho}.

The paper is organized as follows: in Section \ref{Four.Main} the method of Fourier decoupling of \cite{Sar,Bat.Sar.Yom2} is presented in some
detail, next we define the metric span $\omega$ and give our main results. In Section \ref{examples} one specific example is considered in detail, 
illustrating, in particular, the importance of the frequency bound in the general results of Section \ref{Four.Main}. In Section \ref{Exp.Pol} we 
study uniqueness and robustness of non-uniform sampling of exponential polynomials. Finally, in Section \ref{Num.Sim} some results of numerical 
simulations are presented.

\section{Robustness of Fourier Decoupling}\label{Four.Main}
\setcounter{equation}{0}

We consider signals of the form (\ref{equation_decoupling_model}):

$$
F(x)=\sum_{j=1}^k \sum_{q=1}^{q_j} a_{jq} f_j(x-x_{jq}), \ a_{jq}\in \mathbb{C}, \ x_{jq} \in {\mathbb R}.
$$
Here $f_j$ are known, while $a_{jq}, \ x_{jq}$ are the unknown signal parameters, which we want to find from Fourier samples ${\cal F}(F)(s)$ of 
$F$ at certain sample points $s\in {\mathbb R}$. Let ${\cal F}(f_j)$ be the (known) Fourier transforms of $f_j$.

For $F$ of the form (\ref{equation_decoupling_model}) and for any $s \in {\mathbb R}$ we have for the sample of the Fourier transform ${\cal F}(F)$ 
at $s$

\be\label{sample}
{\cal F}(F)(s)=\sum_{j=1}^k \sum_{q=1}^{q_j} a_{jq} e^{-2\pi isx_{jq}}{\cal F}(f_j)(s).
\ee
In the case $k=1$ we could divide the equation (\ref{sample}) by ${\cal F}(f_1)(s)$ and obtain directly a Prony-like equation. However, for $k\geq 2$ 
this transformation usually is not applicable. Instead, in \cite{Bat.Sar.Yom2} we ``decouple'' equations (\ref{sample}) with respect 
to the signals $f_1,\ldots,f_k$ using the freedom in the choice of the sample set $S$. Let \[Z_\l=\bigl\{x\in {\mathbb R}, \ {\cal F}(f_\l)(x)=0\bigr\}\] 
denote the set of zeroes of the Fourier transform ${\cal F}(f_\l)$. For each $j=1,\dots,k$ we take the sampling set $S_j$ to be a subset of the set 
\[W_j=W_j(f_1,\ldots,f_k)=\biggl(\bigcap_{\l\ne j} Z_\l\biggr)\setminus Z_j\] of common zeroes of the Fourier transforms ${\cal F}(f_\l), \ \l\ne j$, but not of 
${\cal F}(f_j)$. For such $S_j$ all the summands in (\ref{sample}) vanish, besides those with the index $j$. Hence we obtain:

\bp\label{proposition_coupling} (\cite{Bat.Sar.Yom2})
Let for each $j=1,\dots,k$ the sampling set $S_j$ satisfy \[S_j=\{s_{j1},\dots,s_{jm_j}\}\subset W_j.\] Then for each
$j$ the corresponding system of equations (\ref{sample}) on the sample set $S_j$ takes the form

\be\label{equation_decoupled1}
\sum_{q=1}^{q_j} a_{jq} e^{-2\pi ix_{jq}s_{j\l}} = c_{j\l}, \quad \l=1,2,\dots, \ s_{j\l} \in S_j,
\ee
where $c_{j\l}=c_{j\l}(F)= {{{\cal F}(F)(s_{j\l})}/{{\cal F}(f_j)(s_{j\l})}}$.
\ep 
These decoupled systems are exactly the same as the fitting systems for exponential polynomials $H_j(s)=\sum_{q=1}^{q_j} a_{jq} e^{-2\pi ix_{jq}s}$ 
on sets $S_j$. So various methods of exponential fitting can be applied (see, for example, \cite{ixaru_exponential_2004, pereyra_exponential_2010,peter2011nonlinear,stoica2005spectral} and references therein.) As it was 
mentioned above, the main problem is that the sample sets $W_j$ may be non-uniform, and/or not sufficiently dense to provide a robust fitting. Indeed,
the zeroes sets $Z_\l$ of the Fourier transforms ${\cal F}(f_\l)$ may be any closed subsets $G_\l$ of $\mathbb R$: it is enough to take ${\cal F}_\l$
to be smooth rapidly decreasing functions on $\mathbb R$ with zeroes exactly on $G_\l$, and to define the signals $f_\l$ as the inverse Fourier 
transforms of ${\cal F}_\l$. In particular, as a typical situation, $Z_\l$ may be arbitrary finite sets or discrete sequences of real points.

\smallskip

The main results of this paper provide a simple ``density'' condition on the sets $W_j(f_1,\ldots,f_k)$ as above, which guarantees a robust least 
square reconstruction of the signal $F$ as in (\ref{equation_decoupling_model}). We need some definitions: 

\smallskip

Let $S$ be a bounded subset of $\mathbb R$, and let $I = [0,R(S)]$ be the minimal interval containing $S$. Let $\lambda \in  \mathbb{R}_+$ be fixed. 
We put $M=M(N,\lambda,R(S)) = N^2-1+\lfloor {{\lambda R(S)}\over {\pi}} \rfloor$, where for a real $A$, \ $\lfloor A \rfloor$  denotes the integer 
part of $A$.

\bd \label{m.span}
For $N\in {\mathbb N}, \ \lambda \in {\mathbb R}_+$, the $(N,\lambda)$-metric span of $S$ is defined as
$$
\omega_{N,\lambda}(S) = \max \{0, \ \sup_{\e > 0} \e [M(\e, S) - M(N,\lambda,R(S))]\},
$$
where $M(\e,S)$ is the $\e$-covering number of $S$, i.e. the minimal number of $\e$-intervals covering $S\cap I$. 
\ed
\bd \label{m.freq}
For each $j=1,\ldots,k$ the maximal frequency $\eta_j$ of the $j$-th equation in the decoupled system (\ref{sample}) is defined by 
$$
\eta_j=\max_{q=1,\ldots,q_j} 2\pi |x_{jq}|.
$$ 
The minimal gap $\sigma_j$ of the $j$-th equation in (\ref{sample}) is defined by 
$$
\sigma_j=\min_{1\leq p<q \leq q_j} 2\pi |x_{jq}-x_{jp}|.
$$
\ed
Now let an interval $I_j=[0,R_j]$ be fixed for each $j=1,\ldots,k$, such that $R_j$ be a point in $W_j$. We take the sampling sets $S_j$ of the form 
$S_j=W_j\cap I_j,$ so $I_j$ is the minimal interval of the form $[0,R]$ containing $S_j$, and $R(S_j)=R_j$. In this paper we shall consider only such 
sampling sets $S_j$. This restriction is not essential, but it significantly simplifies the presentation.

\bd \label{min.div}
For each $j=1,\ldots,k$ the minimal divisor $\kappa_j=\kappa_j(S_j)$ of the $j$-th equation in the decoupled system (\ref{sample}) on $S_j$ is 
defined by
$$
\kappa_j=\min_{s\in S_j} |{\cal F}(f_j)(s)|.
$$
The sample gap $\rho_j$ of the $j$-th equation in (\ref{sample}) on $I_j$ is defined by
$$
\rho_j={{3R_j \sigma_j}\over {2\pi q_j^2(q_j+1)}} \ for \ \eta_j R_j \leq \pi q_j, \ and \ \rho_j={{2\sigma_j}\over {\eta_j q_j(q_j+1)}} \ otherwise.
$$
\ed
\bt \label{sign.main}
Assume that for each $j=1,\ldots,k$ we have $\omega_j:=\omega_{2q_j,\eta_j}(S_j,R_j)>0.$ Then the parameters 
$a_{jq}, \ x_{jq}, \ q=1,\ldots,q_j, \ j=1,\ldots,k$ of the signal $F$ as in (\ref{equation_decoupling_model}) can be uniquely reconstructed via  
the least square solution of the equations (\ref{equation_decoupled1}) on the sample sets $S_j$, assuming that the measured samples of ${\cal F}(F)(s)$  
at all the sample points are exact.

In the case of noisy measurements, with the maximal error of the sample ${\cal F}(F)(s_{j\l})$ for $s_{j\l} \in S_j$ being at most $\delta_j$ 
(sufficiently small), we have the following bounds for the reconstruction errors $\Delta a_{jq}, \ \Delta x_{jq}$ of $a_{jq}, \ x_{jq}$:

\be\label{amppl}
\Delta a_{jq} \leq {2\over {\kappa_j}}\cdot \left(\frac{632 R_j}{\rho_j\omega_j}\right)^{2q_j} \cdot \delta_j,
\ee
\be\label{freqqu}
\Delta x_{jq} \leq {{2}\over {|a_{jq}|\kappa_j}}\cdot \left(\frac{632 R_j}{\rho_j\omega_j}\right)^{2q_j} \cdot \delta_j.
\ee
\et
\pr
This theorem follows directly from Theorem \ref{non.unif.rec} below, which estimates the accuracy of the least square sampling of exponential
polynomials with purely imaginary exponents on a given sampling set $S$. The only adaptation we have to make is that the right hand sides 
$c_{j\l}$ of the equations (\ref{sample}) are given by $c_{j\l}=c_{j\l}(F)= {{{\cal F}(F)(s_{j\l})}/{{\cal F}(f_j)(s_{j\l})}}$, and hence the Fourier 
sampling error is magnified by ${1\over {{\cal F}(f_j)(s_{j\l})}}$. Consequently, the minimal divisor $\kappa_j=\kappa_j(S_j)$ of the $j$-th equation 
in (\ref{equation_decoupled1}) on $S_j$ appears in the denominator of (\ref{amppl}) and (\ref{freqqu}). $\square$

As a corollary we show that if the Fourier zeroes sets $W_j$ are ``sufficiently dense'' then the decoupling approach provides a robust reconstruction
of the signal $F$. The notion of density we introduce below is a very restricted one. Much more accurate definition, involving not only the asymptotic
behavior of $W_j\cap [0,R]$ as $R$ tends to infinity, but also its finite geometry, can be given. We plan to present these results separately. Notice
also a direct connection with the classical Sampling Theory, in particular, with Beurling theorems of \cite{Beu,landau_necessary_1967}. See also \cite{Marvasti,Ole.Ula1} and
references therein.

\bd \label{dens}
Let $S$ be a discrete subset of $\mathbb R_+$. The ``central density'' $D(S)$ of $S$ is defined as 
$D(S)=\limsup_{R\rightarrow \infty} {{|S\cap [0,R]|}\over R}.$ 
\ed
\bc \label{density}
If for $j=1,\ldots,k$ we have $D(W_j) > {\eta_j\over \pi}$ then the decoupling procedure on appropriate sampling sets $S_j\subset W_j$ provides a 
robust reconstruction of the signal $F$.
\ec
\pr
If $D(W_j) > {\eta_j\over\pi}$ then for arbitrarily big $R_j$ we have $|W_j\cap [0,R_j]| > M(2q_j,\eta_j,R_j),$ so taking sufficiently small $\e>0$ we 
conclude that the span $\omega_{2q_j,\eta_j}(S_j,R_j)$ is strictly positive. Application of Theorem \ref{sign.main} completes the proof. $\square$

\section{An example}\label{examples}
\setcounter{equation}{0}

Some examples of Fourier decoupling have been presented in \cite{Sar,Bat.Sar.Yom2}. In the present paper we consider one of these examples
in more detail, stressing the question of robust solvability of the resulting decoupled systems. As everywhere in this paper, we restrict ourselves 
to the case of one-dimensional signals. Some initial examples in dimension two can be found in \cite{Sar,Bat.Sar.Yom2}.

\smallskip

Let $f_1$ be the characteristic function of the interval $[-1,1],$ while we take $f_2(x)=\delta(x-1)+\delta(x+1).$ So we consider signals of the form

\be\label{equation_decoupling_model1}
F(x)= \sum_{q=1}^{N} [a_{1q} f_1(x-x_{1q})+a_{2q} f_2(x-x_{2q})].
\ee
We allow here the same number $N$ of shifts for each of the two signals $f_1$ and $f_2$. Easy computations show that

\[
{\cal F}(f_1)(s)=\sqrt{\frac{2}{\pi}}\frac{\sin s}{s}
\]
and
\[
{\cal F}(f_2)(s)=\sqrt\frac{2}{\pi}\cos s.
\]
So the zeros of the Fourier transform of $f_1$ are the points $\pi n,\ n\in {\mathbb Z}\setminus \{0\}$ and those of $f_2$ are the points 
$({1\over 2}+n)\pi,\ n\in {\mathbb Z}$. These sets do not intersect, so we have $W_1=\{\pi n\}$, and $W_2=\{({1\over 2}+n)\pi\}$, and we can take as 
$S_1, S_2$ any appropriate subsets of these sets. Notice that the central density of $W_1,W_2,$ according to Definition \ref{dens}, is $1\over \pi$. 
By Corollary \ref{density}, if both the maximal frequencies $\eta_1,\eta_2$ are strictly smaller than $1$ then the decoupling procedure on appropriate 
sampling sets $S_j\subset W_j, \ j=1,2,$ provides a robust reconstruction of the signal $F$. Let us show that this condition on the frequencies 
$\eta_1,\eta_2$ is sharp. 

The decoupled systems, given by (\ref{equation_decoupled1}) above, take the form

\be\label{equation_decoupled3}
\sum_{q=1}^{N} a_{1q} e^{-2({1\over 2}+n)\pi^2 ix_{1q}}=\sum_{q=1}^{N}\alpha_qe^{i\phi_qn} = c_{1n}, \ n\in {\mathbb Z},
\ee

\be\label{equation_decoupled4}
\sum_{q=1}^{N} a_{2q} e^{-2n\pi^2 ix_{2q}} = \sum_{q=1}^{N} \beta_q e^{i\psi_qn} = c_{2n}, \ n\in {\mathbb Z},
\ee
where

$\alpha_q=a_{1q}e^{-i\pi^2x_{1q}}, \ \phi_q = -2\pi^2x_{1q}, \ \beta_q=a_{2q}, \ \psi_q = -2\pi^2x_{2q},$

$c_{1n}= {{{\cal F}(F)(({1\over 2}+n)\pi)}/{{\cal F}(f_1)(({1\over 2}+n)\pi)}}$, $c_{2n}= {{{\cal F}(F)(n\pi)}/{{\cal F}(f_2)(n\pi)}}$.

Now we put in the equations (\ref{equation_decoupled3}), (\ref{equation_decoupled4})

$\alpha_1={1\over {2i}}, \alpha_2=-{1\over {2i}}, \alpha_q=0$ for $q=3,\ldots,N$, $\phi_1=\pi, \ \phi_2 = -\pi$,

$\beta_1={1\over {2i}}, \beta_2=-{1\over {2i}}, \beta_q=0$ for $q=3,\ldots,N$, $\psi_1=\pi, \ \psi_2 = -\pi$.

This corresponds to the following shifts and amplitudes in the signal $F$:

$x_{11}= - {1\over {2\pi}}, \ x_{12}= {1\over {2\pi}}, \ a_{11}=-{i\over 2} \ a_{12}={i\over 2}$, and thus $\eta_1=1$.

$x_{21}= - {1\over {2\pi}}, \ x_{22}= {1\over {2\pi}}, \ a_{21}=-{1\over 2} \ a_{22}={1\over 2}$, \ \ $\eta_2=1$.

For this specific signal $F$ the exponential polynomials (\ref{equation_decoupled3}), (\ref{equation_decoupled4}) are both equal to $\sin(\pi s)$, so 
they both vanish identically at all the sampling points $s=n \in {\mathbb Z}$. Thus, allowing $\eta_1=\eta_2=1$ we cannot reconstruct uniquely our 
signals from the samples on the sets $W_1,W_2,$ no matter how many sampling points we take.

\smallskip

On the other hand, put $\eta_1=\eta_2=\eta < 1$. Let us take as $S^m_1$ (respectively, $S^m_2$) the set of points of the form $({1\over 2}+n)\pi$ 
(respectively, $n\pi$) for $n=0,\ldots,m$. We have $R(S^m_1)=({1\over 2}+m)\pi, \ R(S^m_2)=m\pi$. The number of the sample points in each case is 
$m+1.$ So in computing $\omega_{2N,\eta,m\pi}(S^m_2)$ we have $M = 4N^2-1+\lfloor {{\eta m\pi}\over {\pi}} \rfloor \sim 4N^2+\eta m-1$. So we have
$\omega_{N,\lambda,m\pi}(S^m_2)\sim \sup_{\e>0}\e[M(\e,S^m_2)-4N^2+1-\eta m].$ Substituting here $\e < \pi$ tending to $\pi$, we get $M(\e,S^m_2)=m+1$, 
so $\omega(S^m_2) \sim [(1-\eta)m - N^2+2]\pi.$ Essentially the same expression we get for $\omega(S^m_1)$. So the metric spans of $S^m_1$ and 
$S^m_2$ are positive for $m > {{4N^2-1}\over {1 - \eta}}.$ Applying Theorem \ref{sign.main} we conclude that for such $m$ the least square sampling 
on the sets $S^m_1,S^m_2,$ is well posed, and get explicit estimates for its accuracy. It would be very desirable to check the sharpness of this 
conclusion. Our numerical simulations, presented below, provide an initial step in this direction. 

\section{Sampling of Exponential Polynomials}\label{Exp.Pol}
\setcounter{equation}{0}

The decoupling method of \cite{Sar,Bat.Sar.Yom2}, presented in Section \ref{Four.Main} above, reduces the Fourier reconstruction problem 
for signals of the form (\ref{equation_decoupling_model}) to a system of decoupled equations (\ref{equation_decoupled1}), which are, for
each $j=1,\ldots,k,$ the sampling equations for exponential polynomials of the form $H_j(s)=\sum_{q=1}^{q_j} a_{jq} e^{-2\pi ix_{jq}s}$
on sampling sets $S_j$. So from now on we deal with sampling of exponential polynomials, not returning any more to the original 
problem of the Fourier reconstruction of linear combinations of shifts of several signals.  

\subsection{Problem definition and main assumptions}

We study robustness of sampling of exponential polynomials on the real line. Let

\be\label{Basic.Exp.Pol}
H(s)=\sum_{j=1}^{N} a_j e^{\lambda_js}, \ a_j,\lambda_j\in \C, \ s\in {\mathbb R},
\ee
be an exponential polynomial of degree $N$. We consider the following problem.

\smallskip

{\it Given a sampling set $S\subset {\mathbb R},$ can an exponential polynomial $H$ of degree $N$ be reconstructed (i.e. its coefficients 
$a_j,\;\lambda_j$ be recovered) from its known values on $S$? If so, how robust can this reconstruction procedure be with respect to noise 
in the data?}

Let us now elaborate some assumptions we keep below.
\begin{enumerate}
\item In this paper we deal with the case of only purely imaginary exponents $\lambda_j = \imath \phi_j, \ \phi_j \in {\mathbb R}$. This 
assumption, which is satisfied in the case of Fourier reconstruction of the linear combinations of shifts of several signals, i.e. for the 
fitting problem (\ref{equation_decoupled1}) above, strongly simplifies the presentation. We plan to describe the general case of arbitrary 
complex exponents separately.
\item We restrict ourselves to the least square reconstruction method, and do not consider other possible reconstruction schemes.
\item In the noisy setting, we investigate the case of sufficiently small noise levels (for a more detailed explanation of this assumption 
see Theorem \ref{thm:siam} below, and \cite{Bat2,Bat.Yom,Yom,Bat.Yom2}).
\end{enumerate}

As it was shown above, in order to ensure well-posedness of the (even noiseless) reconstruction problem, a certain ``density'' of the sampling 
set $S$ with respect to the frequency set $\{\phi_1,\dots,\phi_N\}$ must be assumed. Accordingly, we assume an explicit upper bound $\lambda$ on 
the frequencies $\phi_j$ and incorporate this bound into the definition of the metric span of the sampling sets (compare Definition \ref{m.span}
above). We shall also assume a lower bound on the minimal distance between the frequencies: $|\phi_j-\phi_i| \geq \Delta.$ Without this assumption 
we cannot bound the accuracy of the reconstruction of the amplitudes $a_j$: indeed, as the exponents $\lambda_j$ of the exponential polynomial 
$H(s)$ as in (\ref{Basic.Exp.Pol}) collide, while the amplitudes $a_j$ tend to infinity in a pattern of divided finite differences, $H(s)$ remains 
bounded on any finite interval (see \cite{Yom,Bat.Yom2}). Accordingly, we shall always assume that for certain fixed $\lambda > 0, \Delta >0$ we
have

\be\label{eq:assumptions-freq}
\max \{\phi_1,\dots,\phi_N\} \leq \lambda, \ \min_{i<j} |\phi_j-\phi_i| \geq \Delta.
\ee
The inequalities \eqref{eq:assumptions-freq} will serve also as the constraints in our least square fitting procedure.

\subsection{Reconstruction by least squares}
Let there be given the sampling set $S=\{s_1,\dots,s_n\}$ of size $n$ and the noisy samples of some unknown exponential polynomial $H$ of degree $N$:
\[
h_k = H(s_k) + \delta_k,\qquad k=1,\dots,n.
\]
According to our assumptions, the noise satisfies
\[
|\delta_k|<\delta,
\]
where $\delta$ is assumed to be sufficiently small.
Let the exponential polynomial
\[
\tilde H(s)= \sum_{j=1}^N \tilde a_j e^{\tilde \lambda_j s},
\]
with
$\tilde a_j \in \mathbb{C}$, $\tilde \lambda_j=i\tilde \phi_j, \ \tilde \phi_j \in {\mathbb R},$ provide the least square fitting of the samples $h_k$, 
under the constraints \eqref{eq:assumptions-freq}. That is,
\[
\left(\tilde a_j,\;\tilde \phi_j\right)=\arg \min_{|\tilde \phi_j|\leqslant \lambda,\;|\tilde \phi_i-\tilde \phi_j| \geq \Delta}
\sum_{k=1}^n \bigl|\sum_{j=1}^N \tilde a_j e^{\tilde \lambda_j s_k}-h_k\bigr|^2
\]

At this stage we do not assume that $\tilde H(s)$ is uniquely defined by the sampling data. Our goal is to estimate the deviations $|a_j-\tilde a_j|$ 
and $|\phi_j-\tilde \phi_j|$ as function of $\Delta,N,n,\lambda,S$ and $\delta$. The approach is as follows:

\begin{enumerate}
\item First we estimate the difference $|H-\tilde{H}|$ at every point $s\in S$, via a simple comparison of the least square deviations for $H$ and $\tilde H$.
\item Then we estimate $|H-\tilde{H}|$ on a certain \emph{interval} $I,$ with $S\subset I$, using discrete version of Turan-Nazarov inequality. At this stage
a major role is played by the metric span of $S$.
\item Now we choose inside the interval $I$ a certain arithmetic progression of points $\bar S=\{s_0,2s_0,\ldots,(2N-1)s_0\}$. The reconstruction problem 
on $\bar S$ is reduced to the standard Prony system. The right hand side of this Prony system, i.e. the values of $\tilde H$ on $\bar S$, would deviate from 
those of $H$ not more than allowed by the estimate of the previous step. Then, the deviations of the reconstructed parameters $\tilde a_j,\;\tilde \phi_j$ from 
the original ones can finally be estimated by the Lipschitz constant of the inverse Prony mapping, as presented in \cite{Bat.Yom} (see Theorem \ref{thm:siam} 
below). An appropriate choice of $s_0$ is possible if we assume (as we do) that the exponents $\phi_j$ do not collide.
\end{enumerate}

The rest of this section is organized as follows. The discrete Tur\'an-Nazarov inequality is presented in Subsection \ref{sub:TN}. The stability estimates for 
the inverse Prony mapping are reproduced in Subsection \ref{sub:Prony}. The formulation of the final estimate and its proof using the above steps are presented 
in Subsection \ref{sub:exp-main-proof}.

\subsection{The discrete Turan-Nazarov inequality}\label{sub:TN}

Let $I=[0,R(S)]$ be the minimal interval containing $S$. Let $N \in \mathbb{N}$ and $\lambda \in  \mathbb{R}_+$ be fixed. We recall that the metric span 
$\omega_{N,\lambda}(S)$ was defined as $\max \{0, \ \sup_{\e > 0} \e [M(\e, S) - M(N,\lambda,R(S)]\},$ where 
$M(N,\lambda,R) = N^2-1+\lfloor {{\lambda R}\over {\pi}} \rfloor$. In this section we shall prove the following special case of the main result of \cite{Fri.Yom}:

\bt \label{d.turan.naz}
Let $ H(s) = \sum_{j=1}^N  a_j e^{\lambda_j s}$ be an exponential polynomial, where
$ a_j \in \C, \ \lambda_j=i \phi_j , \ \phi_j \in {\mathbb R}, \ \ \lambda = \max_{j=1,\ldots,N} | \phi_j| $. Let $S,I$ be as above, with 
$\omega_{N,\lambda}(S)>0$. Then we have
\be \label{tn}
\sup_{I} |H(s)| \le \left( \frac{316 R(S)}{\omega_{N, \lambda}(S) }\right)^{N-1} \cdot \sup_{S} | H(s)|.
\ee
\et

We follow the lines of proof of Theorem 1.3 of \cite{Fri.Yom}. We shall use the following two results from \cite{Naz}.

\bl[Lemma 1.3 in \cite{Naz}, Langer's lemma] \label{lem:langer}
Let $p(z)=\sum_{k=1}^n c_k e^{i\lambda_k z}$ be an exponential polynomial ($0<\lambda_1 < \dots < \lambda_n=\lambda$) not vanishing identically. Then the number of complex zeros of $p(z)$ in an open vertical strip $x_0 < \Re z < x_0 + \Delta$ of width $\Delta$ does not exceed
\[
n-1 + \frac{\lambda \Delta}{2\pi}.
\]
\el

\bt[Theorem 1.5 in \cite{Naz}, the T\'{u}ran's lemma]\label{tur-naz}
Let $p(t)=\sum_{k=1}^n c_k e^{i\lambda_k t}$, where $c_k\in\mathbb{C}$ and $\lambda_1 < \dots < \lambda_n \in \mathbb{R}$. If $E$ is a measurable subset of an interval $I=[0,R]$ then
\be\label{nazarov-ineq}
\sup_{t\in I} \left| p(t) \right| \leqslant \biggl\{ \frac{316R}{\mu\left( E \right)} \biggr\}^{n-1} \sup_{t\in E} \left| p(t) \right|,
\ee
where $\mu$ is the Lebesgue measure on $\mathbb{R}$.
\et

{\bf Proof of Theorem \ref{d.turan.naz}: }
Let $\rho:=\sup_S |H(s)|$ and consider the sublevel set
\[
V_\rho := \{t\in I: |H(s)| < \rho \}.
\]
Define $p(z)=H^2(z)-\rho^2$. It is an exponential polynomial of degree at most $N^2$ with purely imaginary exponents whose absolute values are bounded from above by $2\lambda$. By  Lemma \ref{lem:langer} above the number of solutions to the equation $p(t)=0$ in $I$ (which is equivalent to $|H(s)|=\rho$) is at most
\[
N^2-1+\biggl\lfloor \frac{\lambda R(S)}{\pi} \biggr\rfloor =M(N,\lambda,R(S)).
\]

Therefore, the set $V_\rho$ consists of at most $M(N,\lambda,R(S))$ subintervals $\Delta_i$. Now fix $\epsilon>0$ and consider the $\epsilon$-covering number $M(\epsilon,V_\rho)$. In order to cover each of the $\Delta_i$'s, we need at most ${\mu(\Delta_i)\over\epsilon}+1$ $\epsilon$-intervals. Overall, we get
\be\label{mevrho}
M(\epsilon,V_\rho) \leqslant M(N,\lambda,R(S)) + {\mu(V_\rho)\over\epsilon}.
\ee
By definition of $\rho$, we obviously have that $S\subseteq V_\rho$, therefore $M(\epsilon,S)\leqslant M(\epsilon,V_\rho)$. Substituting this into  \eqref{mevrho}, 
multiplying by $\epsilon$ and taking supremum w.r.t. $\epsilon$ we obtain $\mu(V_\rho) \geqslant \omega_{N,\lambda}(S)$.  Now we just apply Theorem \ref{tur-naz} with $p=H$ and $E=V_\rho$.
$\square$

 Notice that the result of Theorem \ref{d.turan.naz} does not depend at all on the 
minimal distance between the exponents of $ H$, which is crucial in the rest of our estimates, and does not imply any bound on the amplitudes $ a_j$ of 
$H$.  

\subsection{Robustness estimates of inverse Prony mapping}\label{sub:Prony}
Let $x_1,\dots,x_N$ be pairwise distinct complex numbers, and let $a_1,\dots,a_N$ be nonzero complex numbers. In \cite{Bat.Yom} we introduced the 
``Prony map'', ${\cal P}: \mathbb{C}^{2N}\to \mathbb{C}^{2N}$, defined by
\[
{\cal P}(x_1,\dots,x_N, a_1, \dots, a_N) = (m_0, \dots, m_{2N-1}),\qquad m_k = \sum_{j=1}^N a_j x_j^k.
\]
This mapping can be considered as the sampling operator for the exponential polynomial $H(s)=\sum_{j=1}^N a_j x_j^s$ on the integer points 
$s\in \{0,1,\dots,2N-1\}$. In \cite{Bat.Yom} we provided local perturbation estimates for $\cal P$, as follows.

\bt\label{thm:siam}
Let $x_1,\dots,x_N$ be pairwise distinct complex numbers, and let $a_1,\dots,a_N$ be nonzero complex numbers. Let $\vec{x}=(m_0,\dots,m_{2N-1})$ be the 
image of the point $(x_1,\dots,x_N,a_1,\dots,a_N)$ under the Prony map $\cal P$. Let $\delta>0$ be sufficiently small, so that the inverse map 
${\cal P}^{-1}$ is defined in the $\delta$-neighborhood $U$ of $\vec{x}$. Let $\vec{\tilde x}$ be some point in this neighborhood:
\[
\vec{\tilde x}=(m_0+\delta_0,\dots,m_{2N-1}+\delta_{2N-1}),\qquad |\delta_i|<\delta.
\]
 Then the image of $\vec{\tilde x}$ under ${\cal P}^{-1}$ satisfies
 \begin{align}
 \begin{split}
 |a_j-\tilde a_j| &\leqslant C(x_1,\dots,x_N) \delta, \\
 |x_j-\tilde x_j| &\leqslant C(x_1,\dots,x_N) |a_j|^{-1} \delta,
\end{split}
 \end{align}
 where $C(x_1,\dots,x_N)$ depends only on the configuration of the nodes $x_1,\dots,x_N$.
\et
In fact, as we show in \cite{Bat2}, in the case that $x_1,\dots,x_N$ belong to the unit circle, the constant $C$ can be bounded from above by
\be\label{eq:c-bound}
C \leqslant 2\cdot \biggl( \frac{2}{\Lambda} \biggr)^{2N},
\ee
where $\Lambda = \min_{i<j} |x_i-x_j|$. While it is not known if the bound \eqref{eq:c-bound} is sharp, it appears to be reasonably accurate in asymptotic terms, as demonstrated in \cite{donoho}.

 As for the the size $\delta$ of the neighborhood $U$ of the point $\vec{x}=(m_0,\dots,m_{2N-1}),$ where the
inverse map ${\cal P}^{-1}$ is defined, its explicit determination is not straightforward, since the geometry of the Prony map, as well as its
singularities, are rather complicated. In \cite{Yom,Bat.Yom2} we have started algebraic-geometric investigation of the Prony map, and the results
there provide some explicit information on $\delta$.

\subsection{Accuracy of least squares sampling}\label{sub:exp-main-proof}

The following is our main result on the least square sampling of $H$ on $S$:

\bt \label{non.unif.rec}
Let $H(s)=\sum_{j=1}^N a_j e^{\imath\phi_j s}$ be an a-priori unknown exponential polynomial satisfying 
$\max |\phi_j|\leq \lambda, \ \min |\phi_i-\phi_j|\geq \Delta$ for some fixed $\lambda,\Delta$. Let there be given the noisy samples $h(s)$ of $H(s)$ 
on a finite set $S=\{s_1,\dots,s_n\}\subset \mathbb{R}$, with the noise bounded by $\delta,$ i.e. $|h(s_\l)-H(s_\l)|\leq \delta, \ \l=1,\ldots, n.$ 
Assume that $\omega(S):=\omega_{2N,\lambda}(S)>0$. Then, for sufficiently small $\delta$, the amplitudes $\tilde a_j$ and the frequencies $\tilde \phi_j$ 
of the least square fitting exponential polynomial $\tilde H(s)$ satisfy, for $j=1,\ldots,N,$ the following inequalities:
\be\label{ampl}
|a_j-\tilde a_j| \leq 2\sqrt{2n}\cdot \left( \frac{632 R(S)}{\rho \omega(S) }\right)^{2N} \cdot \delta,
\ee
\be\label{frequ}
|\phi_j-\tilde \phi_j| \leq 2\sqrt{2n} \ |a_j|^{-1} \cdot \left( \frac{632 R(S)}{\rho \omega(S) }\right)^{2N} \cdot \delta, 
\ee
where $\rho={{3R(S) \Delta}\over {2\pi N^2(N+1)}}$ for $\lambda R(S) \leq \pi N$, and $\rho={{2\Delta}\over {\lambda N(N+1)}}$ otherwise. In 
particular, in the case of zero noise, the least square reconstruction of $\tilde H(s)$ on $S$, under the constraints as above, is unique, up to a 
transposition of the indices.
\et
\pr
First of all, let us establish the following easy bound:

\bl\label{quadr.dev}
For $s\in S$ we have $|\tilde H(s) - H(s)| \ \leq \sqrt{2n} \delta$.
\el
\pr
Indeed, the quadratic deviation $\sigma(H,h)$ of $H$ from $h$ on $S$ does not exceed $n\delta^2$, where $n$, as above, denotes the number of elements
in $S$. Since $\tilde H(s)$ is the exponential polynomial of the least square deviation from $h$, we have $\sigma(\tilde H,h)\leq n\delta^2,$ which
directly implies $\sigma(\tilde H,H)\leq 2n\delta^2$ and hence $|\tilde H(s) - H(s)| \ \leq \sqrt{2n} \delta,$ for each $s\in S$. $\square$

Now we get directly the following bound:

\bc\label{cor.bound}
For $H,S$ and $I=[0,R(S)]$ the minimal interval containing $S$ we have

\be
\sup_{I} |\tilde H(s) - H(s)| \ \leq \left( \frac{316 R(S)}{\omega_{2N,\lambda}(S) }\right)^{2N-1} \cdot \sqrt{2 n} \delta.
\ee
\ec
\pr
We notice that by Lemma \ref{quadr.dev} we have $\sup_{S} |\tilde H(s) - H(s)| \ \leq \sqrt{2 n} \delta$. Substituting into Theorem \ref{d.turan.naz}
(which is applied to the exponential polynomial $H(s) - \tilde H(s)$ of degree $2N$ with purely imaginary exponents, bounded in absolute value by
$\lambda$), we get the required bound. $\square$

The bound of Corollary \ref{cor.bound} {\it does not imply by itself any bound on the amplitudes $a_j$}. They may tend
to infinity, as the exponents collide, following the pattern of divided finite differences (see \cite{Yom,Bat.Yom2}). So the continuation of the
proof incorporates the a priori known lower bound $\Delta$ on the differences between the exponents of $H$. We get estimates of the reconstruction
accuracy of $a_j$ and $\lambda_j$ via solving an appropriate auxiliary Prony system, and applying Theorem \ref{thm:siam} above.

\smallskip

Let $I=[0,R(S)]$ be as above. Fix certain $s_0 \in (0, {R\over {2N}}]$ and consider the points $s_0,2s_0,\ldots,(2N)s_0 \in I$. We denote 
$\nu_k=H(ks_0), \ k=0,1,\ldots,$ the values of $H$ at the points $ks_0$. We get
\be\label{Aux.Prony1}
H(ks_0)=\sum_{j=1}^N a_j e^{\l_j ks_0}=\sum_{j=1}^N a_j x_j^k= \nu_k, \ k \in {\mathbb Z},
\ee
where $x_j=e^{\lambda_js_0}=e^{i\phi_js_0}.$ So for each choice of $s_0 \in (0, {R\over{2N}}]$ we obtain a Prony system

\be\label{Aux.Prony2}
\sum_{j=1}^N a_j x_j^k= \nu_k, \ k=0,\ldots,2N-1,
\ee
which is satisfied by $a_j$ and $x_j=e^{i\phi_js_0}, \ j=1,\ldots,N$. It is well known that if $x_i\ne x_j$ for $i\ne j,$ then the solution $a_j,$
$x_j=e^{i\phi_js_0}, \ j=1,\ldots,N$ of (\ref{Aux.Prony2}) is unique, up to a permutation of the indices. Moreover,
the robustness of the solutions of (\ref{Aux.Prony2}) with respect to the perturbations of the right-hand side, is determined by the mutual
distances $|x_i- x_j|, \ i\ne j$ (see \cite{Bat.Yom, Bat.Yom1} and Subsection \ref{sub:Prony}). So our next goal is to choose $s_0 \in (0, {R\over {2N}}]$ 
in such a way that $\Lambda=\min_{i\ne j} |x_i- x_j|$ be sufficiently large. To achieve this goal we have to find $s_0$ such that all the angles 
$\Delta_{i,j}s_0$ are separated from the integer multiples $2\pi m, \ m\in {\mathbb Z}$, where $\Delta_{i,j}=|\phi_j-\phi_j|.$ 

An easy example shows that there are ``bad'' choices of $s_0$: assume that the frequencies $\phi_j$ in $H$ are of the form $\phi_j= s\cdot 2\pi m_j,$ with
$s\in {\mathbb R}, \ m_j \in {\mathbb Z}, \ m_i\ne m_j$ for $i \ne j$. Then for $s_0={1\over s}$ we have $x_1=x_2=\ldots = x_N.$ The next lemma shows that 
most choices of $s_0$ are good, assuming that $\Delta=\min_{i<j} \Delta_{i,j}$ is not zero.

\bl\label{arithmetics}
Let $\bar R, \ q_1,\ldots q_r \in {\mathbb R}_+$ be given, with $q=\min q_\l, \ Q=\max q_\l$ There exists $s_0 \in (0, \bar R]$ such 
that all the angles $q_\l s_0, \ \l=1\ldots,r,$ are separated from the integer multiples $2\pi m, \ m\in {\mathbb Z}$ by at least $\bar h$, defined as 
$\bar h={{\bar R q}\over {4r}}\leq {\pi\over 4}$ for $Q\bar R \leq \pi$, and as $\bar h={{\pi q}\over {3rQ}}\leq {\pi\over 3}$ for $Q\bar R > \pi$.
\el
\pr By assumptions we have $q_\l \leq Q$. Hence for each $\l=1,\ldots,r$ the interval $q_l\cdot(0, \bar R]$ contains at most ${{Q \bar R}\over {2\pi}}+1$ 
integer multiples $2\pi m$. For $h>0$ let $U(h)$ denote the $h$-neighborhood of these points. Denote by $\mu_1$ the standard Lebesgue measure on $\mathbb{R}$. 
We have $\mu_1(U(h))\leq h[{{Q \bar R}\over {\pi}}+2]$. Now let $V_\l(h)$ denote the set of those $s\in (0, \bar R]$ for which $q_\l s \in U(h)$. We conclude 
that $\mu_1(V_\l(h)(h))\leq {h\over {q_\l}}[{{Q \bar R}\over {\pi}}+2]\leq {h\over {q}}[{{Q \bar R}\over {\pi}}+2]$. Finally, denoting $V(h)$ the set of 
the points $s\in (0, \bar R]$ for which $q_\l s \in U(h)$ for at least one index $\l=1,\ldots,r$, we get 
$\mu_1(V(h))\leq {{rh}\over {q}}[{{Q \bar R}\over {\pi}}+2]$. If for some $h$ we have $\mu_1(V(h))< |(0, \bar R]|= \bar R$, then there exists 
$s_0 \in (0, \bar R]$ such that all the angles $q_l s_0, \ \l=1,\ldots,r,$ are separated from the integer multiples $2\pi m, \ m\in {\mathbb Z}$ at least 
by $h$.

Now we consider two cases: $Q\bar R \leq \pi$ and $Q\bar R > \pi$. In the first case $\mu_1(V(h))\leq {{3rh}\over {q}},$ and
the inequality $\mu_1(V(h))\leq \bar R$ is valid with $h=\bar h={{\bar R q}\over {4r}}\leq {\pi\over 4}.$ In the second case
$\mu_1(V(h))\leq {{3rh}\over {q}}{{Q \bar R}\over {\pi}},$ and the inequality $\mu_1(V(h))\leq \bar R$ is valid, starting with 
$h=\bar h = {{\pi q}\over {3rQ}}\leq {\pi\over 3}.$ This completes the proof of the lemma. $\square$

In our case of the angles $\Delta_{i,j}$ and $s_0 \in (0, {R\over {2N}}]$ we have, respectively, 
$\bar R = {R\over {2N}}, \ r={{N(N+1)}\over 2}, Q=2\lambda, \ q=\Delta$. Applying Lemma \ref{arithmetics} we obtain the following result:

\bc \label{distances}
There exists $s_0 \in (0, {R\over {2N}}]$ such that all the angles $\Delta_{i,j}\cdot s_0, \ 1<i<j\leq N$ are separated from the integer multiples 
$2\pi m, \ m\in {\mathbb Z}$ by at least $\bar h$, defined as $\bar h={{R \Delta}\over {2N^2(N+1)}}\leq {\pi\over 4}$ for $\lambda R \leq \pi N$, and 
as $\bar h={{2\pi \Delta}\over {3\lambda N(N+1)}}\leq {\pi\over 3}$ for $\lambda R > \pi N$. Accordingly, the minimal distance 
$\min |x_i-x_j|, \ i\ne j,$ between the points $x_j=e^{i\phi_js_0}, j=1,\ldots,N$ in (\ref{Aux.Prony2}) is at least $\rho={3\over \pi} \bar h,$ which
is $\rho={{3R \Delta}\over {2\pi N^2(N+1)}}$ for $\lambda R \leq \pi N$, and $\rho={{2\Delta}\over {\lambda N(N+1)}}$ for $\lambda R > \pi N$.
\ec
\pr The result on the separation of the angles follows directly from Lemma \ref{arithmetics}. The result for the distances follows from the fact that
always $\bar h \leq {\pi \over 3}$. $\square$ 

\smallskip

Now we can complete the proof of Theorem \ref{non.unif.rec}. We fix $s_0$ whose existence is guaranteed by Lemma \ref{arithmetics}, and form the
Prony system, which is satisfied by the parameters of $H$:

\be\label{Aux.Prony3}
\sum_{j=1}^N a_j x_j^k= \nu_k, \ k=0,\ldots,2N-1, \ x_j=e^{\lambda_j s_0,}
\ee
with $\nu_k=H(ks_0)$ the values of $H$ at the points $ks_0, \ k=0,\ldots,2N-1$. Notice that these values are not exactly known. However, by
Corollary \ref{cor.bound} we know that

\be\label{next}
\sup_{J} |\tilde H(s) - H(s)| \ \leq \left( \frac{316 R}{\omega_{2N,\lambda}(S) }\right)^{2N-1} \cdot \sqrt{2 n} \delta,
\ee
where $\tilde H(s)=\sum_{j=1}^N \tilde a_j e^{\tilde \lambda_j s}$ is the  polynomial of the least square approximation on $S$. In particular,
denoting $\tilde \nu_k= \tilde H(ks_0), \ k=0,\ldots,2N-1,$ the values of $\tilde H$ at the points $ks_0,$ we get
$|\tilde \nu_k - \nu_k| \leq \left( \frac{316 R}{\omega_{2N,\lambda}(S) }\right)^{2N-1} \cdot \sqrt {2n} \delta.$

Now the parameters $\tilde a_j, \tilde \lambda_j$ of $\tilde H$ satisfy the Prony system

\be\label{Aux.Prony4}
\sum_{j=1}^N \tilde a_j \tilde x_j^k= \tilde \nu_k, \ k=0,\ldots,2N-1, \ \tilde x_j=e^{\tilde \lambda_j s_0.}
\ee
Finally we apply Theorem \ref{thm:siam}  to Prony system (\ref{Aux.Prony3}) and its perturbation (\ref{Aux.Prony4}), taking into account the
expression (\ref{eq:c-bound}) for the constant $C$ in Theorem \ref{thm:siam}. Noticing that the distances between the nodes $x_j$ of the 
unperturbed Prony system (\ref{Aux.Prony3}) are bounded from below by $\rho$ via Corollary \ref{distances}, we arrive at \eqref{ampl} and 
\eqref{frequ}. Uniqueness of reconstruction for $\delta = 0$ follows directly from \eqref{ampl} and \eqref{frequ}. This completes the proof 
of Theorem \ref{non.unif.rec}. $\square$

\subsection{Estimating $\omega_{N,\lambda}(S)$: some examples}\label{actual.calc}

The metric span $\omega_{N,\lambda}(S)$ can be explicitly computed in many important cases. In particular, we have the following simple result:

\bp \label{span.comp}
Let $N,\lambda$ be fixed. Assume that a subset $S \subset {\mathbb R}$ with $R(S)=R$ contains $M(N,\lambda,R)+1$ points, and let $\eta$ be the
minimal distance between the neighboring points in $S$. Then $\o_{N,\lambda}(S)= \eta.$
\ep
\pr
For $\e\geq \eta$ we have $M(\e, S) - M(N,\lambda,R) \leq 0$. For $\e < \eta$ this difference is $1$. Hence the supremum in Definition
\ref{m.span} is achieved as $\e$ tends to $\rho$ from the left. $\square$

\bc \label{extr.span.comp}
Let $N,\lambda$ be fixed. Assume that a subset $S \subset {\mathbb R}$ contains $M(N,\lambda,R)+1$ points. Then
$\o_{N,\lambda}(S)\leq {R\over {M(N,\lambda,R)}},$ and this value is achieved only for $S$ consisting of $M(N,\lambda,R)+1$ points at the
distance ${{R}\over {M(N,\lambda,R)}}$ one from another.
\ec
\pr
For the equidistant configuration the minimal distance $\eta$ between the neighboring points in $S$ is ${R\over {M(N,\lambda,R)}}$. Otherwise
$\eta$ is strictly smaller. $\square$.

\medskip

Now let us consider equidistant configurations with a larger number of sampling points.

\bp \label{equidist.span.comp}
Let $N,\lambda$ be fixed. Assume that a subset $S \subset {\mathbb R}$ contains $m+1\geq M(N,\lambda,R)+1$ points at the distance
${R\over m}$ from one another. Then
\be\label{many.pts}
\o_{N,\lambda}(S)= \biggl({R\over m}\biggr)[m+1-M(N,\lambda,R)].
\ee
\ep
\pr
For the equidistant configuration $S$ the minimal distance $\eta$ between the neighboring points in $S$ is ${R\over {m}}$. On the other
hand, for each $\e< \eta$ we have $M(\e,S)=m+1,$ while for $k\eta \leq \e \leq (k+1)\eta, \ k=1,2,...,$ we have $M(\e,S)={{m+1}\over k}.$
An easy computation then shows that the supremum of $\e[M(\e,S)- M(N,\lambda,R)]$ is achieved for $\e$ tending to $\eta$ from the left, and
it is equal to $({R\over {m}})[m+1-M(N,\lambda,R)].$ $\square$

\smallskip

\begin{remark}\label{rem:npoints}
As substituted into the expression of Theorem \ref{non.unif.rec}, the results above imply the corresponding bound for the accuracy of the
least square reconstruction on $S$. In particular, the expression \eqref{many.pts} above seems to provide a non-trivial recommendation for the
choice of the number of equidistant sample points inside a given interval $I$. Indeed, for $m=M(N,\lambda,R)$ we get
$\omega(S)={R\over {M(N,\lambda,R)}}$. But for $m=2 M(N,\lambda,R)$ we get 
\[
\omega(S)={R\over {2M(N,\lambda,R)}}\bigl[M(N,\lambda,R)+1\bigr],
\]
which is approximately ${R\over 2}$ for large $M(N,\lambda,R)$ - improvement by $M\over 2$ times. For $m$ tending to infinity $\omega(S)$ tends
to $R$, so we do not achieve any essential improvement any more. Thus the recommendation may be to take $m$ of order $KM_d$ with $K$
between, say $2$ and $5$. 
\end{remark}

\smallskip

\begin{remark}\label{rem:collision}
Combining Theorem \ref{non.unif.rec} and Proposition \ref{span.comp} we can also predict the rate of the degeneration of the reconstruction
problem on $S$ as two points of $S$ collide. By the same method we can analyse also the cases of more complicated collisions between the
sampling points.
\end{remark}

\begin{figure}
\begin{center}
\includegraphics[scale=0.5,trim=2cm 6cm 2cm 7.5cm, clip=true]{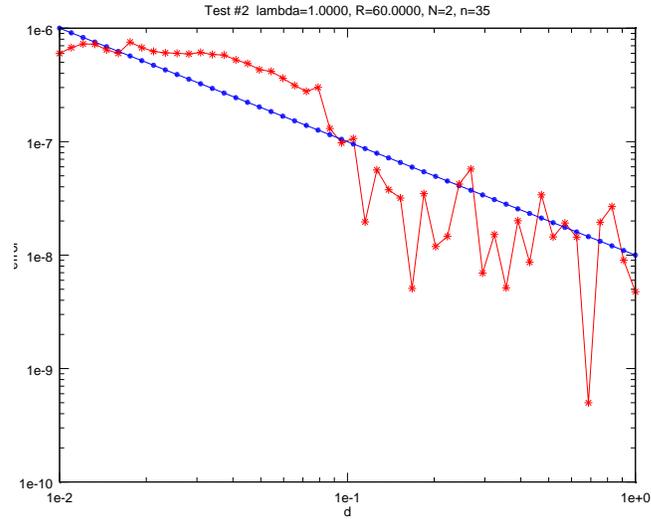}
\end{center}
\caption[Experiment 1]{In this experiment, we changed the mutual distance $d$ between the subsequent points of $S$, while keeping the two endpoints fixed. $\varepsilon_2 = 10^{-5},\;\lambda=1,\;R=60,\;N=2$. The size of $S$ is $n=35$. The error is plotted versus the value of $d$ in red. For comparison, the value ${1\over d}$ is plotted in blue.}\label{fig:exp2}
\end{figure}

\begin{figure}
\begin{center}
\includegraphics[scale=0.5,trim=2cm 6cm 2cm 7cm, clip=true]{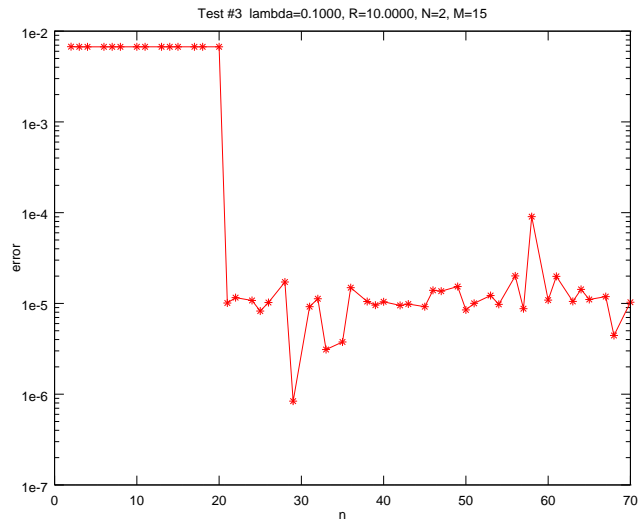}
\end{center}
\caption[Experiment 2]{In this experiment, we increased $n$, the number of points in $S$, keeping the range (i.e. the value of $R$) fixed. $\varepsilon_2 = 10^{-2},\;\lambda=0.1,\;R=10,\;N=2$. Here $M(2N,\lambda,R)=15$. The error is plotted versus the value of $n$.}
\label{fig:exp3}
\end{figure}

\section{Numerical simulations}\label{Num.Sim}
In this section we present results of initial numerical experiments. Our goal in these very preliminary simulations has been to numerically investigate the qualitative dependence of the reconstruction error on the geometry of the sampling set $S$. Our results below are indeed qualitatively consistent with the bounds of Theorem \ref{non.unif.rec}.

In all the experiments presented in Figures \ref{fig:exp2} and \ref{fig:exp3} below, we have fixed an a-priori randomly chosen exponential polynomial $H(s)$, and modified the sampling set $S$ according to the description of each experiment below. The sampling values $\{H(s_i),\; s_i \in S\}$ have been perturbed by the (random) amount $\varepsilon_1 \sim 10^{-8}$. Subsequently, the least-squares approximation to $H(s)$ has been obtained by the standard sequential quadratic programming algorithm (implemented by the function  \verb!sqp! in GNU Octave environment). The initial values for the algorithm have been taken to be equal to the true values perturbed by the (random) amount $\varepsilon_2$, specified in each experiment below. We have plotted the recovery error for one of the frequencies (specifically, $|\Delta \phi_2|$).

In the first experiment we changed the distance $d$ between $s_2,\dots,s_{n-1}$, while keeping the endpoints $s_1,\;s_n$ (and thereby the value of $R$) fixed. The number of points was chosen to be exactly $n=M(2N,\lambda,R)+1$. According to Proposition \ref{span.comp}, in this case we have $\omega(S)=d$. As can be seen in Figure \ref{fig:exp2}, the error is roughly proportional to ${1\over \omega(S)}$.

In the second experiment, we have kept the endpoints of the set $S$ fixed ($0$ and $R$), while increasing the number $n$ of (equispaced) points in $S$. According to Figure \ref{fig:exp3}, a significant improvement in accuracy appears when the number of samples passes $M(2N,\lambda,R)$ which is $15$ in this case.

\bibliographystyle{plain}

\begin{thebibliography}{10}

\bibitem{adk} B.~Adcock, M.~Gataric, and A.C.~Hansen
\newblock{On stable reconstructions from univariate nonuniform Fourier measurements.}
\newblock {\em Preprint. Arxiv: 1310.7820}

\bibitem{Bat} D.~Batenkov.
\newblock {Complete Algebraic Reconstruction of Piecewise-Smooth Functions from Fourier Data}.
\newblock {\em To appear in Mathematics of Computation.}

\bibitem{Bat2} D.~Batenkov.
\newblock {Decimated Generalized Prony systems}.
\newblock {\em Preprint. Arxiv:1308.0753.}

\bibitem{Bat.Sar.Yom1}
D.~Batenkov, N.~Sarig, and Y.~Yomdin.
\newblock {An ``algebraic" reconstruction of piecewise-smooth functions from
integral measurements}.
\newblock {\em Functional Differential Equations}, 19(1-2):9--26, 2012.

\bibitem{Bat.Sar.Yom2}
D.~Batenkov, N.~Sarig, and Y.~Yomdin.
\newblock {Decoupling of Reconstruction Systems for Shifts of Several Signals}.
\newblock {\em Proc. of Sampling Theory and Applications (SAMPTA)}, 2013. 

\bibitem{Bat.Yom1} D.~Batenkov and Y.~Yomdin.
\newblock {Algebraic reconstruction of piecewise-smooth functions from Fourier data}.
\newblock {\em Proc. of Sampling Theory and Applications (SAMPTA)}, 2011.

\bibitem{Bat.Yom} D.~Batenkov and Y.~Yomdin.
\newblock {On the accuracy of solving confluent Prony systems}.
\newblock {\em SIAM J. Appl. Math.}, 73(1):134--154, 2013.

\bibitem{Bat.Yom2} D.~Batenkov and Y.~Yomdin.
\newblock {Geometry and Singularities of the Prony Mapping}.
\newblock {\em To appear in Proceedings of 12th International Workshop on Real and Complex Singularities}, 2013.

\bibitem{Beu} A. Beurling. 
\newblock {Balayage of Fourier-Stiltjes Transforms}. 
\newblock {\em The collected Works of Arne Beurling, Vol.2, Harmonic Analysis}.
Birkhauser, Boston, 1989.

\bibitem{donoho} D.L.~Donoho.
\newblock {Superresolution via sparsity constraints}.
\newblock {\em SIAM Journal on Mathematical Analysis,}
\newblock{23(5):1309--1331, 1992.}

\bibitem{Vet} P.L. ~Dragotti, M. ~Vetterli and T. ~Blu.
\newblock {Sampling Moments and Reconstructing Signals of Finite Rate of Innovation:
Shannon Meets Strang-Fix},
\newblock {\em IEEE Transactions on Signal Processing,}
\newblock {
Vol. 55, Nr. 5, Part 1, pp. 1741-1757, 2007.}

\bibitem{Fri.Yom}
O.~Friedland and Y.~Yomdin.
\newblock { An observation on Tur\'an-Nazarov inequality}.
\newblock {\em Studia Mathamatica},
\newblock{218(1), pp. 27--39, 2013}.
\newblock {\em DOI 10.4064/sm218-1-2}

\bibitem{gedalyahu2011multichannel}
K.~Gedalyahu, R.~Tur, and Y.C. Eldar.
\newblock Multichannel sampling of pulse streams at the rate of innovation.
\newblock {\em IEEE Transactions on Signal Processing}, 59(4):1491--1504, 2011.

\bibitem{gust.Moments}
B.~Gustafsson, C.~He, P.~Milanfar and M.~Putinar.
\newblock Reconstructing planar domains from their moments.
\newblock{\em Inverse Problems}, 16(4):1053--1040, 2000.

\bibitem{ixaru_exponential_2004}
Liviu~Gr Ixaru and Guido~Vanden Berghe.
\newblock {\em Exponential Fitting}.
\newblock Springer, May 2004.

\bibitem{landau_necessary_1967}
H.~Landau.
\newblock Necessary density conditions for sampling and interpolation of
  certain entire functions.
\newblock {\em Acta Mathematica}, 117(1):37--52, 1967.

\bibitem{Marvasti} F. ~Marvasti.
\newblock{\em Nonuniform sampling: theory and practice}.
\newblock{Springer, 2001.}

\bibitem{Naz} F.L. Nazarov.
\newblock {Local estimates of exponential polynomials and their applications to
  inequalities of uncertainty principle type}.
\newblock {\em St Petersburg Mathematical Journal}, 5(4):663--718, 1994.

\bibitem{Ole.Ula1} A. Olevski, A. Ulanovski.
\newblock {Near critical density irregular sampling in Bernstein spaces}.
\newblock {\em Mathematisches Forschungsinstitut Oberwolfach gGmbH}, Oberwolfach 
Preprints (OWP) 2013-16, ISSN 1864-7596.

\bibitem{pereyra_exponential_2010}
Victor Pereyra and Godela Scherer.
\newblock {\em Exponential data fitting and its applications}.
\newblock Bentham Science Publishers, January 2010.

\bibitem{peter2011nonlinear}
T.~Peter, D.~Potts, and M.~Tasche.
\newblock Nonlinear approximation by sums of exponentials and translates.
\newblock {\em SIAM Journal on Scientific Computing}, 33(4):1920, 2011.

\bibitem{rao1992mbp}
B.D. Rao and K.S. Arun.
\newblock {Model based processing of signals: A state space approach}.
\newblock {\em Proceedings of the IEEE}, 80(2):283--309, 1992.

\bibitem{sig.ack}
N.~Sarig and Y.~Yomdin.
\newblock {Signal Acquisition from Measurements via Non-Linear Models}.
\newblock {\em Mathematical Reports of the Academy of Science of the Royal
  Society of Canada}, 29(4):97--114, 2008.

\bibitem{Sar}
N.~Sarig.
\newblock {\em {Algebraic reconstruction of "shift-generated" signals from
  integral measurements}}.
\newblock PhD thesis, {Weizmann Institute of Science}, 2010.

\bibitem{stoica2005spectral}
P.~Stoica and R.L.~Moses.
\newblock {\em Spectral analysis of signals}.
\newblock Pearson/Prentice Hall, 2005.

\bibitem{Yom} Y.~Yomdin.
\newblock Singularities in algebraic data acquisition.
\newblock {\em Real and Complex Singularities (M. Manoel, {MC} Romero Fuster,
  {CTC} Wall, eds.), London Mathematical Society Lecture Notes}, 380:378--396,
  2010.


\end{thebibliography}

\end{document}